\documentclass[a4paper,11pt]{article}

\usepackage{amssymb,latexsym,amsmath,amsthm,amsfonts,graphics}

\def\prf{\vskip -0.1cm\noindent{\bf Proof\quad }}
\def\prfend{\hfill{$\Box$}\vskip 0.2cm}

\def\D{\{0,1\}}

\def\N{\mathbb{N}}
\def\A{\mathbb A}

\def\1{\mathbf 1}

\def\F{\mathcal{F}}

\def\NN{\mathcal{N}}

\def\R{\mathbb{R}}
\def\Z{\mathbb{Z}}

\def\P{\mathcal{P}}

\newtheorem{lem}{Lemma}

\newtheorem{thm}{Theorem}

\newtheorem{dfn}{Definition}
\newtheorem{exa}{Example}

\numberwithin{equation}{section}

\begin{document}
\parindent=1em

\title{A characterization of Ramsey property with application to 
dynamical systems\\
{\small Dedicated to Prof. Benjamin Weiss for his 85th birthday}}
\author{Teturo KAMAE\\
\small{Satakedai 5-9-6, 565-0855 Japan}\\
\small{kamae@apost.plala.or.jp}}
\date{}
\maketitle

\date{}
\maketitle

\begin{abstract}
The infinitary Ramsey theorem says that in any finite partition 
$\{V^i;~i=1,2,\cdots,n\}$ of the the family of $k$-subsets of 
$\N$, there exists $V^i$ having an infinite set 
$\NN\subset\N$ such that any $k$-set from $\NN$ is in $V^i$. 
This property of $V^i$ is called the Ramsey property. 
We characterize this property in term of ultrafilter.  
\end{abstract}

\section{Introduction}

For a set $S$, denote $\F(S)$ the family of nonempty subsets of $S$, 
and $\F_k(S)$ the family of subsets of $S$ of cardinality $k$ 
($k=1,2,\cdots$). Let $\N=\{0,1,2,\cdots\}$ and  
$\xi\in\F_k(\N)$.  Arrange the elements 
in $\xi$ in the order so that $\xi=\{\xi_0<\xi_1<\cdots<\xi_{k-1}\}$. 
Thus, we identify the set $\xi$ as this with 
$(\xi_0,\xi_1,\cdots,\xi_{k-1})\in\Delta_k$, where 
$$
\Delta_k=\{(x_0,x_1,\cdots,x_{k-1})\in\N^k;~x_0<x_1<\cdots<x_{k-1}\}.
$$
Hence, 
\begin{equation}
\F_k(\N)=\Delta_k.  
\end{equation}

An infinite set $V$ in $\Delta_k$ is said to have the 
{\it Ramsey property} if there exists an infinite increasing 
sequence $\NN=\{x_0<x_1<\cdots\}\subset\N$ such that 
for any increasing sequence $\{i_0<i_1<\cdots<i_{k-1}\}\subset\N$, 
it holds that $(x_{i_0},x_{i_1},\cdots,x_{i_{k-1}})\in V$. 

The following Ramsey Theorem is well known. 
\begin{thm}{\rm \cite{R}}
For any finite partition $\{V_i;~i=1,2,\cdots,n\}$ of $\Delta_k$, 
there exists $V_i$ which has the Ramsey property. 
\end{thm}

It is an interesting problem to characterize the Ramsey property  
in another term. In this paper, we characterize it 
in term of {\it ultrafilters}.  

For any countably infinite set $\Sigma$, let $\beta\Sigma$ be the 
Stone-Cech compactification of $\Sigma$, that is, 
$\beta\Sigma$ is the set of ultrfilters on $\Sigma$. The set of principal ultrafilters is identified with $\Sigma$, so that 
$\beta\Sigma\setminus\Sigma$ is the set of non-principal 
ultrafilters on $\Sigma$. 
In general, for $\chi_1\in\beta\Sigma_1$ and 
$\chi_2\in\beta\Sigma_2$, 
the product $\chi_1\times\chi_2\in\beta(\Sigma_1\times\Sigma_2)$ 
is defined in Definition 1. Note that this definition is not 
symmetric between $\Sigma_1$ and $\Sigma_2$. That is, 
$V\in\chi_1\times\chi_2$ does not implies 
$\tilde{V}\in\chi_2\times\chi_1$, where 
$\tilde{V}=\{(y,x)\in\Sigma_2\times\Sigma_1;~(x,y)\in V\}$. 
Although the associative law holds, that is, 
$$
(\chi_1\times\chi_2)\times\chi_3=
\chi_1\times(\chi_1\times\chi_2).
$$
For $\chi\in\beta\N$, we consider the power 
$$
\chi^k=
\underbrace{\chi\times\chi\times\cdots\times\chi}_k\in\beta\N^k.
$$
 
By the definition, $\Delta_k\in\chi^k~(k=1,2,\cdots)$ 
holds for any $\chi\in\beta\N\setminus\N$. 
We prove in Lemma 7 that any set $V\in\chi^k$ with 
$V\subset\Delta_k$ and $\chi\in\beta\N\setminus\N$ 
has the Ramsey property. This converse also holds. 

For $V\subset\Delta_k$, let $\P_V$ be the set of 
$\lambda\in\beta\Delta_k$ such that $V\in\lambda$. 
Then, any finite partition $\{V^i;~i=1,2,\cdots,n\}$ of $\Delta_k$ 
define a finite partition $\{\P_{V^i};~i=1,2,\cdots,n\}$
of $\beta\Delta_k$. 

Our main result is that 
\begin{thm}{\rm (Main Theorem)} 
An infinite set $V\subset\Delta_k$ has the Ramsey property 
if and only if $\P_V$ contains a $\chi^k$ for some 
$\chi\in\beta\N\setminus\N$. Hence, for any finite partition 
$\{V_i;~i=1,2,\cdots,n\}$ of $\Delta_k$, there exists $V_i$ 
having the Ramsey property. 
\end{thm}

Let us give an easy example. Let $\{V_1,V_2,V_3\}$ be a partition 
of $\Delta_2$ such that 
\begin{align*}
&V_1=\{(x,y)\in\Delta_2;~x,y~\mbox{are odd}\},\\
&V_2=\{(x,y)\in\Delta_2;~x,y~\mbox{are even}\},\\
&V_3=\{(x,y)\in\Delta_2;~x,y~\mbox{have different parities}\}.
\end{align*}
Then $V_1$ and $V_2$ have the Ramsey property, but $V_3$ 
does not. On the other hand, for an odd 
$\chi_1\in\beta\N\setminus\N$ and even 
$\chi_2\in\beta\N\setminus\N$, we have $\chi_1^2\in\P_{V_1}$ 
and $\chi_2^2\in\P_{V_2}$, while $\P_{V_3}$ does not contain 
any of $\chi^2$ for $\chi\in\beta\N\setminus\N$. 

We give a way of giving finite partitions of $\Delta_k$ 
related to dynamical systems. 

Let $\A$ be a set with $2\le\#\A<\infty$. Let 
$\Omega\subset\A^\N$ be a nonempty closed set. For 
$\omega\in\A^\N$ and 
$\xi=(\xi_0,\xi_1,\cdots,\xi_{k-1})\in\Delta^k$, define 
$$\omega[\xi]:=(\omega(\xi_0),\omega(\xi_1),\cdots,\omega(\xi_{k-1}))
\in\A^k.$$ 
Using this, we define 
$$
\Omega[\xi]:=\{\omega[\xi];~\omega\in\Omega\}\in\F(\A^k). 
$$
Consider the partition of $\Delta_k$ by the value of 
$\Omega[\xi]$ for $\xi\in\Delta_k$. 
That is, 
$$
\{\{\xi\in\Delta_k;~\Omega[\xi]=\Lambda\}\ne\emptyset;~
\Lambda\subset\A^k\}.
$$
This partition is called the $\Omega$-$k$-{\it partition}. 

The $\Omega$-$k$-partitions are specially interesting in the case 
that $\Omega$ is a shift invariant symbolic representation of 
a geometrical dynamical system. In this case, the Ramsey property 
corresponds to {\it infinitesimal moves at infinity} of the 
geometrical dynamical system as discussed in T.Kamae \cite{K4}. 

\section{Product of a non-nonprincipal ultrafilter}

The product of a non-principal ultrafilter (\cite{K4}) plays an 
important role in this paper, so let us review it.  

\begin{dfn}{\rm 
Let $\Sigma_i$ be infinite sets for $i=1,2$ and $\chi_i\in\beta\Sigma_i~(i=1,2)$. We define 
$\chi_1\times\chi_2\in\beta(\Sigma_1\times\Sigma_2)$ as
$$\chi_1\times\chi_2=\{U\subset\Sigma_1\times\Sigma_2;~
\{x\in\Sigma_1;~\{y\in\Sigma_2;~(x,y)\in U\}\in\chi_2\}\in\chi_1\}.
$$
}\end{dfn}

Note that $\chi_1\times\chi_2=\chi_1'\times\chi_2'$ implies  
$\chi_1=\chi_1'$ and $\chi_2=\chi_2'$. 
Also, note that this product operation is not symmetric between the first entry and the second entry. That is, $U\in\chi_1\times\chi_2$ 
and $\tilde{U}\in\chi_2\times\chi_1$ is not equivalent, where 
$\tilde{U}=\{(y,x);~(x,y)\in U\}$. In fact, 
$\{(x,y)\in\N^2;~x<y\}\in\chi\times\chi$ but $\{(y,x)\in\N^2;~x<y\}\notin\chi\times\chi$ for any $\chi\in\beta\N\setminus\N$, since any cofinte set is in $\chi$, but any finite set is not in $\chi$. 
Although the following associative law follows clearly from 
Definition 1.

\begin{lem} 
For any $\chi_i\in\beta\Sigma_i~(i=1,2,3)$, 
It holds that 
$$(\chi_1\times\chi_2)\times\chi_3=\chi_1\times(\chi_2\times\chi_3)
\in\beta(\Sigma_1\times\Sigma_2\times\Sigma_3).$$
Hence, the power $\chi^k\in\beta(\Sigma^k)~(k=2,3,\cdots)$ 
for $\chi\in\beta\Sigma$ is defined without ambiguity. 
\end{lem}

\begin{dfn}{\rm
Let $k\ge2$ and 
$S=\{s_0<s_1<\cdots<s_{l-1}\}\subset\{0,1,\cdots,k-1\}$.  
The projection $\pi_{k,S}:\N^k\to\N^l$ is defined as  
$$
\pi_{k,S}(x_0,x_1,\cdots,x_{k-1})
=(x_{s_0},x_{s_1},\cdots,x_{s_{l-1}})$$ 
for any $(x_0,x_1,\cdots,x_{k-1})\in\N^k$. 

We simply denote 
$$\pi^i:=\pi_{k,\{0,1,\cdots,i-1\}}~\mbox{and}~
\pi_i:=\pi_{k,\{k-i,k-i+1,\cdots,k-1\}}$$
for any $1\le i\le k$. 
}\end{dfn}

\begin{dfn}{\rm
Let $1\le l<k$ and $V\subset\N^k$. 
For $\xi=(\xi_0,\xi_1,\cdots,\xi_{l-1})\in\pi^l(V)$, define  
$$
\Pi_\xi(V)=\{(x_0,x_1,\cdots,x_{k-l-1})\in\N^{k-l};~
(\xi_0,\cdots,\xi_{l-1},x_0,\cdots,x_{k-l-1})\in V\}.
$$
}\end{dfn}

The following 2 Lemmas are easy from Definition 1 
(see also \cite{K4}). 

\begin{lem}
For any $\chi\in\beta\N\setminus\N$, $k\ge1$ and $V_0\in\chi$, 
we have $V_0^k\cap\Delta_k\in\chi^k$. 
\end{lem}

\begin{lem}
Let $\chi\in\beta\N\setminus\N$ and $V\in\chi^k$, where $k\ge 2$. 
Let $1\le l<k$ and $S$ be as in {\rm Definition 2}. Then, we have 
$\pi_{k,S}V\in\chi^l$. Also for any $U\in\chi^l$, we have 
$\pi_{k,S}^{-1}U\in\chi^k$.
\end{lem}

\begin{dfn}{\rm 
Let $k\ge2$. We call $V\subset\N^k$ {\it decreasing} if 
the following 2 conditions hold:\\
$(1)$ $\pi^{k-1}V\supset\pi_{k-1}V$.\\
$(2)$ For any 
$$\xi=(\xi_0,\xi_1,\cdots,\xi_{k-2})\in\pi^{k-1}(V)~\mbox{and}~ 
\eta=(\eta_0,\eta_1,\cdots,\eta_{k-2})\in\pi^{k-1}(V)$$
with $\xi\le\eta$, that is, 
$\xi_0\le\eta_0,~\xi_1\le\eta_1,~\cdots,\xi_{k-2}\le\eta_{k-2}$, 
it holds that $\Pi_\xi(V)\supset\Pi_\eta(V)$.  
}\end{dfn}

\begin{dfn}{\rm 
For any $\chi\in\beta\N\setminus\N$ and $2\le k<\infty$, we call 
$V\in\chi^k$ with $V\subset\Delta_k$ $\chi$-{\it pure} if 
for any $\xi\in\pi^l(V)$ with $1\le l<k$, we have    
$\Pi_\xi(V)\in\chi^{k-l}$.     
}\end{dfn}

\begin{lem}
Let $\chi\in\beta\N\setminus\N$ and $2\le k<\infty$. 
Then for any $V\in\chi^k$, there exists 
$W\subset V$ such that $W$ is $\chi$-pure.  
\end{lem}
\prf 
We may and do assume that $V\in\Delta_k$. 
By Definition 1, 
\begin{equation}
W_0:=\{\xi_0\in\pi^1(V);~\Pi_{\xi_0}(V)\in\chi^{k-1}\}\in\chi.
\end{equation}
Let $W'=\cup_{\xi_0\in W_0}\{\xi_0\}\times\Pi_{\xi_0}(V)$. 
Then, $W'\subset V$ since $W_0\subset\pi^1(V)$. Moreover, 
$W'\in\chi^k$. If $k=2$, this implies that $W'$ is $\chi$-pure. 

We use the induction on $k\ge3$. Assume that the statement holds 
for $k-1$. Then we can take $\Gamma_{\xi_0}\subset\Pi_{\xi_0}(V)$ 
in (2.1) such that $\Gamma_{\xi_0}\in\chi^{k-1}$ and is $\chi$-pure.  
Hence $W:=\cup_{\xi_0\in W_0}\{\xi_0\}\times\Gamma_{\xi_0}$ 
satisfies that $W\subset W'\subset V$ and $W\in\chi^k$ by 
Definition 1. Moreover. $W$ is $\chi$-pure by the construction. 
\prfend

\begin{lem}
Let $\chi\in\beta\N\setminus\N$ and $k\ge2$.  
Let $V^{i}\in\chi^k~(i=1,2,\cdots,m)$ be $\chi$-pure. 
Then for any $l=1,2,\cdots,k-1$, we have 
$$\pi^l\left(\bigcap_{i=1,2,\cdots,m}V^i\right)
=\bigcap_{i=1,2,\cdots,m}\pi^l(V^i).$$
\end{lem}
\prf
In general, it holds that 
$$\pi^l\left(\bigcap_{i=1,2,\cdots,m}V^i\right)
\subset\bigcap_{i=1,2,\cdots,m}\pi^l(V^i).$$
Let $\xi\in\bigcap_{i=1,2,\cdots,m}\pi^l(V^i)$. 
We prove that it belongs to the left hand side in the above. 
Since 
$V^i~(l=1,2,\cdots,k-1)$ are $\chi$-pure, it holds that 
$\Pi_\xi(V^i)\in\chi^{k-l}$ and hence 
$\bigcap_{i=1,\cdots,l}\Pi_\xi(V_i)\ne\emptyset$. 
Then, $\xi\eta\in V^{i}$ holds 
for any $i=1,2,\cdots,l$ by taking 
$\eta\in\bigcap_{i=1,\cdots,l}\Pi_\xi(V_i)$. Thus, 
$\xi\in\pi^l\left(\bigcap_{i=1,2,\cdots,m}V^i\right)$.  
\prfend

\begin{lem}
Let $\chi\in\beta\N\setminus\N$ and $k\ge2$.  
Then for any $V\in\chi^k$, there exists $W\subset V$ with 
$W\in\chi^k$ such that $W$ is decreasing.  
\end{lem}
\prf
By Lemma 2 and 4, we may assume that $V\in\Delta_k$ and is 
$\chi$-pure. We define $W'\subset V$ such that $W'$ is $\chi$-pure 
and satisfies the condition (1) in Definition 4 as follows: 
$$
W'=\bigcap_{i=0}^{k-1}\N^i\times\pi^{k-i}V. 
$$
This is because by Lemma 5, 
\begin{multline*}
\pi^{k-1}W'=\bigcap_{i=0}^{k-1}\pi^{k-1}(\N^i\times\pi^{k-i}V)\\
=\N^{k-1}\cap\bigcap_{i=0}^{k-2}\pi^{k-1}(\N^i\times\pi^{k-i}V)
=\bigcap_{i=0}^{k-2}(\N^i\times\pi^{k-i-1}V)
\end{multline*}
while
\begin{multline*}
\pi_{k-1}W'\subset\bigcap_{i=0}^{k-1}\pi_{k-1}(\N^i\times\pi^{k-i}V)
=\pi_{k-1}V\cap\bigcap_{i=1}^{k-1}(\N^{i-1}\times\pi^{k-i}V)
\\=\pi_{k-1}V\cap\bigcap_{i=0}^{k-2}(\N^i\times\pi^{k-i-1}V)
\subset\bigcap_{i=0}^{k-2}(\N^i\times\pi^{k-i-1}V)=\pi^{k-1}W'. 
\end{multline*} 

Now we define $W\subset W'$ which satisfies the condition (2)  
as well as (1). Let 
$$
W=\bigcup_{\xi\in\pi^{k-1}W'}
\left(\{\xi\}\times
\bigcap_{{\eta\in\pi^{k-1}W'}\atop{\eta\le\xi}}\Pi_\eta(W')\right).
$$
Then, $W\in\chi^k$ since $\Pi_\eta(W')\in\chi$ for any 
$\eta\in\pi^{k-1}W'$ with $\eta\le\xi$, and there are finitely many 
$\eta$ as this. It is clear that $W$ is decreasing.  
\prfend

\begin{lem}
Let $\chi\in\beta\N\setminus\N$ and $V\in\chi^k$ with 
$k\ge2$ and $V\subset\Delta_k$. 
Then, there exists a sequence of nonnegative integers 
$a_0<a_1<a_2<\cdots$ such that 
$(a_{i_0},a_{i_1},\cdots,a_{i_{k-1}})\in V$ for any 
sequence of nonnegative integers $i_0<i_1<\cdots<i_{k-1}$. 
\end{lem}
\prf
By Lemma 4 and 6, we may assume that $V$ is pure and decreasing. 
We construct $a_0<a_1<a_2<\cdots$ as required as follows. 
Take $\xi^{(0)}\in\pi^{k-1}(V)$ arbitrary. Then, 
take $a_0\in\Pi_{\xi^{(0)}}(V)$ arbitrary. 
Let $\xi^{(1)}=\pi_{k-1}(\xi^{(0)}a_0)$. 
Now for $n=1,2,\cdots$, 
take $a_n\in\Pi_{\xi^{(n)}}(V)$ arbitrary and 
let $\xi^{(n+1)}=\pi_{k-1}(\xi^{(n)}a_n)$.
In this way, we construct $a_0<a_1<a_2<\cdots$. 

Since $\xi^{(n)}a_n\in V$ and 
$\xi^{(n)}=\pi_{k-1}(\xi^{(n-1)}a_{n-1})$, we have 
\begin{align*}
&\xi^{(n)}a_n=\pi_{k-1}(\xi^{(n-1)}a_{n-1})a_n
=\pi_{k-1}(\pi_{k-1}(\xi^{(n-2)}a_{n-2})a_{n-1})a_n\\
=&\pi_{k-2}(\xi^{(n-2)}a_{n-2})a_{n-1}a_n=\cdots
=\pi_0(\xi^{(n-k)}a_{n-k})a_{n-k+1}\cdots a_{n-1}a_n\\
=&a_{n-k+1}\cdots a_{n-1}a_n\in V 
\end{align*}
for any $n\ge k-1$.   
Now take an arbitrary subsequence 
$(a_{i_0},a_{i_1},\cdots,a_{i_{k-1}})\in\Delta_k$ with 
$0\le i_0<i_1<\cdots<i_{k-1}$. Denote 
$$(b_0,b_1,\cdots,b_{k-1})=(a_{i_0},a_{i_1},\cdots,a_{i_{k-1}}).$$

Then, $\xi^{(i_0)}b_0=\xi^{(i_0)}a_{i_0}\in V$. Since
$$\xi^{(i_1)}b_1=\xi^{(i_1)}a_{i_1}\in V~\mbox{and}~
\pi_{k-1}(\xi^{(i_0)}b_0)\le\xi^{(i_1)},$$
we have $\pi_{k-1}(\xi^{(i_0)}b_0)b_1\in V$ because of 
the decreasing property of $V$. 
Hence, we have $\pi_{k-2}(\xi^{(i_0)})b_0b_1\in V$. 
In the same way, since
$$\xi^{(i_2)}b_2=\xi^{(i_2)}a_{i_2}\in V~\mbox{and}~
\pi_{k-1}(\pi_{k-2}(\xi^{(i_0)})b_0b_1)\le\xi^{(i_2)},$$
we have $\pi_{k-1}(\pi_{k-2}(\xi^{(i_0)})b_0b_1)b_2\in V$. 
Hence, we have $\pi_{k-3}(\xi^{(i_0)})b_0b_1b_3\in V$. 
Continuing this process, we finally have 
$$b_0b_1\cdots b_{k-1}=\pi_0(\xi^{(i_0)})b_0b_1\cdots b_{k-1}\in V.$$
\prfend

\section{Main Theorem and Examples}

{\bf Main Theorem} (Theorem 2):
An infinite set $V\subset\Delta_k~(k=1,2,\cdots)$ has the 
Ramsey property 
if and only if $\P_V$ contains a $\chi^k$ for some 
$\chi\in\beta\N\setminus\N$. Hence, for any finite partition 
$\{V^i;~i=1,2,\cdots,n\}$ of $\Delta_k$, there exists $V_i$ 
having the Ramsey property. 
\vspace{0.5em}
\prf
The `if' part is exactly Lemma 7. We prove the `only if' part. 

Let $V\subset\Delta_k$ have the Ramsey property.  
Then, there exists an infinite set 
$\NN=\{x_0<x_1<\cdots\}\subset\N$
such that such that 
for any increasing sequence $\{i_0<i_1<\cdots<i_{k-1}\}\subset\N$, 
it holds that $(x_{i_0},x_{i_1},\cdots,x_{i_{k-1}})\in V$. 
Let $\chi$ be any element in $\beta\N\setminus\N$ 
such that $\NN\in\chi$. 
Then, $V\in\chi^k$ holds by Definition 1. 
\prfend

\begin{dfn}{\rm 
Let $\Omega$ be a nonempty closed subset of $\A^\N$ for a 
finite set $\A$ with $\#\A\ge2$. 
Let the $\Omega$-k-partition of $\Delta_k~(k\ge2)$ be
$$
\{\{\xi\in\Delta_k;~\Omega[\xi]=\Lambda\}\ne\emptyset;~
\Lambda\subset\A^k\}.
$$
For any $\chi\in\beta\N\setminus\N$, let $\chi^k\in\P_V$ 
with $V=\{\xi\in\Delta_k;~\Omega[\xi]=\Lambda\}$. 
We define $\Omega[\chi^k]=\Lambda$. 
This coincides with the value of the natural extension of 
$\Omega[\cdot]$ at $\chi^k$ and it is equivalent to 
$\{\xi\in\Delta_k;~\Omega[\xi]=\Lambda\}\in\chi^k$. 
The set of $\Lambda\subset\A^k$ such that 
$\{\xi\in\Delta_k;~\Omega[\xi]=\Lambda\}$ has the Ramsey property 
is denoted by $Sol_k(\Omega)$. 
}\end{dfn} 

\begin{dfn}{\rm 
For $x\in\R$, let 
$\parallel x\parallel:=\min_{n\in\Z}|x-n|$ and for 
$x=(x_0,x_1)\in\R^2$, 
$\parallel x\parallel:=
\sqrt{\parallel x_0\parallel^2+\parallel x_1\parallel^2}$. 
For 3 distinct points $x,y,z\in(\R/\Z)^2$ which are 
close enough each other, we define 
$\mbox{Radius}(x,y,z)\in(0,\infty]$ 
to be the radius of the circle passing these points.  
}\end{dfn} 

The following Lemma follows from the property of any ultrafilter, 
say $\chi$  that 
for any finite partition of $\N$, $\chi$ contains just one of them. 

\begin{lem}{\rm (\cite{K4,K5})} Let $\chi\in\beta\N\setminus\N$. \\
$(1)$ For an irrational $\alpha\in\R$,  
$\beta\in\R/\Z$ is determined so that 
$$\{n\in\N;~\parallel n\alpha-\beta\parallel<\varepsilon\}\in\chi$$ 
for any $\varepsilon>0$. 
Then, we denote $\chi\alpha:=\beta$. 
Moreover in this case, 
\begin{align}
either~~&\{n\in\N;~0<n\alpha-\beta<\varepsilon~(mod~\Z)\}\in\chi\nonumber\\
\mbox{or}~~&\{n\in\N;~-\varepsilon<n\alpha-\beta<0~(mod~\Z)\}\in\chi
\end{align}
holds. Correspondingly, we define more precisely that $\chi\alpha=\beta_+$ 
or $\chi\alpha=\beta_-$. \\
$(2)$ Let $\vec{f}=(f_1,f_2)\in\R^2$ be such that $1,f_1,f_2$ are 
linearly independent over the rational field. 
Then, ${\rm Radius}_{\vec{f}}(\chi)\in(0,\infty)$ is defined so that 
$$
\{(n_0,n_1,n_2)\in\Delta^3;~|{\rm Radius}(n_0\vec{f},n_1\vec{f},n_2\vec{f})
-{\rm Radius}_{\vec{f}}(\chi)|<\varepsilon\}\in\chi^3
$$
for any $\varepsilon>0$. In addition, ${\rm Rad}_{\vec{f}}(\chi)=\infty$ if 
$$\{(n_0,n_1,n_2)\in\Delta_3;~|{\rm Radius}(n_0\vec{f},n_1\vec{f},n_2\vec{f})|
>N\}\in\chi^3$$
for any $N>0$.   
Moreover just like $(3.1)$, we can determine $R_+$ or $R_-$ in the case  
${\rm Radius}_{\vec{f}}(\chi)<\infty$.  
\end{lem}

\begin{exa}{\rm (\cite{K4})
Let $\alpha$ be an irrational number. Let $a,b,c$ be real numbers 
with $a<b<c<a+1$. For $\gamma\in\R/\Z$, define 
$\omega_\gamma\in\{0,1,2\}^\N$ by 
$$\omega_\gamma(n)=\left\{\begin{array}{lll}0&\gamma+n\alpha\in[a,b)&(mod~\Z)\\1&\gamma+n\alpha\in[b,c)&(mod~\Z)\\2&\gamma
+n\alpha\in[c,a+1)&(mod~\Z)\end{array}\right..$$
Let $\Omega$ be the closure of $\{\omega_\gamma;~\gamma\in\R/\Z\}$ 
and k be a positive integer. Then,we have  
$$Sol_k(\Omega)=\{\Omega[\chi^k];~\chi\in\beta\N\setminus\N\}
=\{\Xi^k_+,~\Xi^k_-\},$$
where
\begin{align*}
&\Xi^k_+=\{u^iv^{k-i};~1\le i\le k,~u,v\in\{0,1,2\}\mbox{ are either }v=u\mbox{ or }v=u-1~(mod~3)\}\\
&\Xi^k_-=\{u^iv^{k-i};~1\le i\le k,~u,v\in\{0,1,2\}\mbox{ are either }v=u\mbox{ or }v=u+1~(mod~3)\}.
\end{align*}

For any $\chi\in\beta\N\setminus\N$, let $\chi\alpha=\theta_-$ 
for some $\theta\in\R/\Z$. 
Let 
\begin{align*}
&0<\varepsilon<\delta:=\min\{b-a,c-b,a+1-c\},~\mbox{and}\\
&U=\{(x_0,\cdots,x_{k-1})\in\Delta_k;
~\theta-\varepsilon<x_0\alpha<\cdots<x_{k-1}\alpha<\theta
~(mod~1)\}.
\end{align*}
By Definition 1, $U\in\chi^k$ and there exists $\xi\in U$ such that 
$\Omega[\chi^k]=\Omega[\xi]$. Let $\xi=\xi_0\cdots\xi_{k-1}$. Then since 
$\theta-\delta<\xi_0\alpha<\cdots<\xi_{k-1}\alpha<\theta~(mod~1)$, 
for any $\gamma\in\R/\Z$, $\gamma+\xi_0\alpha,\cdots,
\gamma+\xi_{k-1}\alpha$ are arranged in the order and contained in an open interval of length 
$\varepsilon$. Hence, they stay in one interval of 
$[a,b),[b,c),[c,a+1)$, or at most, move from one interval to the next 
interval. Therefore, $\omega_\gamma[\xi]\in\Xi^k_-$. Moreover, by a choice of $\gamma\in\R/\Z$, $\omega_\gamma(x_0)\cdots\omega_\gamma(x_{k-1})$ can be any element in $\Xi^k_-$. Thus, we have $\Omega[\chi^k]=\Omega[\xi]=\Xi^k_-$. 
In the same way, if $\chi\alpha=\theta_+$ for some $\theta\in\R/\Z$, 
that is, for any $\varepsilon>0$
$$
\{n\in\N;~n\alpha\in[\theta,~\theta+\varepsilon)~(mod~1)\}\in\chi,
$$ 
then we have $\Omega[\chi^k]=\Xi^k_+$. 

Thus, we have $Sol_k(\Omega)=\{\Omega[\chi^k];~\chi\in\beta\N\setminus\N\}=\{\Xi^k_+,\Xi^k_-\}$.  
}\end{exa}

\begin{exa}{\rm (\cite{K4})
Let $\alpha=0110100110010110\cdots\in\D^\N$ be the Thue-Morse sequence. That is, $\alpha(n)=0$ if and only if the number of 0 in the 2-adic representation of $n$ is even. Let $T:\D^\N\to\D^\N$ be the shift and 
$\Omega\subset\D^\N$ be the closure of $\{T^n\alpha;~n\in\N\}$. 
Let $k\ge2$. Then for any $\chi\in\beta\N\setminus\N$, we'll prove 
that $\Omega[\chi^k]=\{\D^k\}$. Thus, the partition having the 
Ramsey property is exactly that corresponding to $\D^k$. 

For the proof refer \cite{K4}. 
}\end{exa}

\begin{exa}{\rm (\cite{K5})
Let $X=(\R/\Z)^2$ and $\vec{f}=(f_0,f_1)\in\R^2$ satisfies that $1,f_0,f_1$ are linearly independent over the rational field. Let 
$$D=\left\{x=(x_0,x_1)\in X;
~\sqrt{\parallel x_0^2\parallel+\parallel x_1^2\parallel}\le1/2\right\}.$$
Define $\omega_0\in \D^\N$ by  
$$\omega_0(n)=\left\{\begin{array}{cc}1&\mbox{if }n\vec{f}\in D~(mod~\Z^2)\\
0&\mbox{otherwise}\end{array}\right.~(n\in\N)$$
and let $\Omega\subset\D^\N$ be the closure of $\{T^n\omega_0;~n\in\N\}$, where $T$ is the shift on $\D^\N$. 

Let $k=3,4,\cdots$. Then we have
$$
Sol_k(\Omega)=\{\Omega[\chi^k];~\chi\in\beta\N\setminus\N\}
=\{\Xi^k,\Lambda^k\},
$$
where 
\begin{align*}
&\Xi^k=\{0^i1^j0^h\in\{0,1\}^k;~i,j,l\in\N,~i+j+l=k\}\\
&\Lambda^k=\Xi^k\cup\{1^i0^j1^l\in\{0,1\}^k;~i,j,l\in\N,~i>0,~l>0,~i+j+l=k\}.
\end{align*}
Moreover, if ${\rm Radius}_{\vec{f}}(\chi)>1/2$ or 
${\rm Radius}_{\vec{f}}(\chi)=(1/2)_+$, then $\Omega[\chi^k]=\Xi^k$ and 
if ${\rm Radius}_{\vec{f}}(\chi)<1/2$ or 
${\rm Radius}_{\vec{f}}(\chi)=(1/2)_-$, then $\Omega[\chi^k]=\Lambda^k$

Let us prove this. Take $\varepsilon>0$ small enough. There exists 
$S\in\chi$ such that $S\vec{f}$ is contained in a $\varepsilon$-neighborhood 
of $\chi\vec{f}$. Moreover, the radius of any 3 points in $S\vec{f}$ is 
larger than $1/2$ in the former case, hence no 3 points $n_0<n_1<n_2$ in 
$S$ satisfy that 
$x+n_0\vec{f}\in D,~x+n_1\vec{f}\notin D,~x+n_2\vec{f}\in D$ since the 
radius of $D$ is $1/2$ . 
Hence, no element of the form $1^i0^j1^l$ with 
$i,j,l\in\N,~i>0,~l>0~i+j+l=k$ is in $\Omega[\chi^k]$. On the other hand, 
by taking the above $x$ adequately, any element in $\Xi$ is realized as an element in $\Omega[\chi^k]$. Thus, $\Omega[\chi^k]=\Xi^k$.
   
In the latter case, in addition to any element in $\Xi^k$ any element in 
$\Lambda$ is also realized as an element in $\Omega[\chi^k]$ since the 
the radius of any 3 points in $S\vec{f}$ is smaller than $1/2$ and 
for any 3 points $n_0<n_1<n_2$ in $S$, 
$x+n_0\vec{f}\in D,~x+n_1\vec{f}\notin D,~x+n_2\vec{f}\in D$ is possible 
with some $x\in(\R/\Z)^2$. Thus, $\Omega[\chi^k]=\Xi^k\cup\Lambda^k$.
}\end{exa}

\hspace{1em}[3][4] can be obtained at 
$\langle$https://teturokamae.web.fc2.com/e-kamae.htm$\rangle$

\end{document}